\documentclass[12pt,twoside]{amsart}
\usepackage{amsmath, amsthm, amscd, amsfonts, amssymb, graphicx, color}
\usepackage[bookmarksnumbered, colorlinks, plainpages]{hyperref}

\newtheorem{theorem}{\bf \textsf Theorem}[section]
\newtheorem{proposition}[theorem]{\bf \textsf {Proposition}}
\newtheorem{lemma}[theorem]{\bf \textsf Lemma}
\newtheorem{corollary}[theorem]{\bf \textsf Corollary}

\numberwithin{equation}{section}

\begin{document}

\begin{center}
{\large \bf  {A new approach  to  Baer and dual Baer modules with some applications}} \\
\end{center}

\begin{center}

\vskip 0.6 true cm

{\bf  N. Ghaedan \footnote{\small
n.ghaedan@math.iut.ac.ir} } and  {\small M.R. Vedadi \footnote{\small
mrvedadi@cc.iut.ac.ir. Corresponding author}\\
Department of Mathematical Sciences, Isfahan University of
Technology, Isfahan, 84156-83111, IRAN.} \\

\end{center}

\vskip 0.6 true cm

\noindent {\bf  {ABSTRACT.}} {\small \  Let $R$ be a ring. It is  proved that an $R$-module $M$ is Baer
(resp. dual Baer) if and only if every exact sequence
$0\rightarrow X\rightarrow M\rightarrow Y\rightarrow 0$ with
$Y\in$ Cog$(M_R)$ (resp. $X\in$ Gen$(M_R)$) splits. This shows that being (dual) Baer is a Morita invariant property. As more applications, the Baer condition for the
$R$-module $M^+ $ = Hom$_{\Bbb Z}(M,{\Bbb Q}/{\Bbb Z})$ is investigated and shown that $R$ is a von Neumann regular ring, if $R^+$ is a Baer $R$-module.
Baer modules with (weak) chain conditions are studied and determined when a  Baer (resp. dual baer) module
 is  a direct sum of mutually orthogonal prime (resp. co-prime) modules.  While   finitely generated dual Baer modules over commutative rings  is shown to be semisimple, finitely generated Baer modules over  commutative domain are studied. In particular, if $R$ is    commutative hereditary Noetherian domain  then a  finitely generated $M_R$  is Baer  if and only if it is projective or semisimple. Over a right duo perfect ring, it is shown that every (dual) Baer modules is    semisimple.
  } \\

\vskip 0.4 true cm

\hrule

\section{\bf Introduction}
\vskip 0.4 true cm

Throughout rings will have  unit elements and modules will be
right unitary. A ring $R$ is said to be {\em Baer} if for  every non-empty subset $X$ of $R$, the right  annihilator $X$ in $R$ is of the form $eR$ for some $e = e^2\in R$.
Baer rings play an important role in the theory of rings of operators in functional analysis; see\cite{kapla1968} and \cite{berb1972}. The concept of Baer ring was extended to modules by S.T. Rizvi and
C. S. Roman in \cite{rizv1}.    A module $M_R$ is called {\em Baer}
if for  every non-empty subset $X$ of  End$_R(M)$, the right  annihilator $X$ in $M$ is a direct summand of $M_R$. Baer modules and their generalizations have been studied among many other works.
Every Baer module $M$ is a D2-module (i.e., if $M/A$ isomorphic to a direct summand of $M$ then $A$ is a direct summand of $M$), see \cite{D3}, \cite{D4}, \cite{direct qbaer2016} and \cite{birken 2020} for recent works on the subjects. The dual notion
of the Baer modules was introduced and studied in \cite{dual.baer2010} where a module $M_R$ is  {\em dual-Baer} if for every $N\leq M_R$, the right ideal
Hom$_R(M,N)$ of End$_R(M)$ is generated by an idempotent element. These modules are known to have the C2-property (i.e., in $M$, every submodule isomorphic to a direct summand is a direct summand), see \cite{C3},\cite{C4},
\cite{cyclicC3}, \cite{w.dual baer2015} and \cite{calci-harmanci 2019} for some recent works on the dual Baer modules  and important generalizations of them. In this paper, we first give  new
characterizations  for a (dual) Baer module of which many are categorical; see Theorem \ref{Baer-exact} and Proposition \ref{preserv direc-rej-1}. These implies that being (dual) Baer is a Morita invariant property and
also co-retractable Baer (resp. retractable dual Baer) modules are semisimple; see Corollary \ref{morita} and Theorem \ref{ret and coret}. Baer (dual Baer) modules with some chain conditions (such as finite uniform dimension) are
studied and determined when a  Baer (resp. dual baer) module
 is  a direct sum of mutually orthogonal prime (resp. co-prime) modules $\bigoplus_{\i\in I} M_i$ such that  Hom$_R(M_i,M_j) = 0$ ($i\neq j$)(Theorems \ref{baer t.dim} and \ref{dual baer t.dim}). While   finitely generated dual Baer modules over commutative rings is shown to be semisimple (Theorem \ref{f.g daul baer}), finitely generated Baer modules over   commutative domain studied in Theorem \ref{comm-domain fg}.
Among other things, it is shown that in  Theorem \ref{perfect- (dual)baer},  over a right duo perfect ring, every (dual) Baer modules is    semisimple. In the last section, we apply our characterization of Baer modules to investigate conditions on $M_R$ under which  the character left $R$-module $M^+$ is Baer. Any unexplained terminology and all the basic results on
rings and modules that are used in the sequel can be found in \cite{AF} and \cite{Lam2}.\\


\section{\bf  Categorical characterization for (Dual) Baer modules}
\vskip 0.4 true cm

If $M_R$ is a nonzero module with End$_R(M) = S$, we shall use the
following notation. If $N\leq M_R$ and $I\leq S_S$, we let $N^c =
$ Hom$_R(M, N)$ and $I^e = $ $IM$. If $K$ and $L$ are two $R$-modules then Tr$(K,L)$ means $\sum\{ f(K) \ | \ f: K\rightarrow L\}$ and Rej$(K,L)$ means
$\bigcap\{ \ker f \ | \ f: K\rightarrow L\}$. If $M_R$ is a module then the class of $R$-modules that are generated (resp. co-generated) by $M$ is denoted by Gen$(M_R)$ (resp. Cog$(M_R)$).
If $_TM_R$ is a bimodule , $A\subseteq T$, $N\subseteq M$ and $B\subseteq R$, then the right annihilator $A$ in $M$ (resp. $N$ in $R$) is canonically defined and denoted by $r_M(A)$
(resp. $r_R(N)$). The left annihilators are similarly denoted by $l_M(B)$ and $r_T(N)$.
 In this section, we consider a characterization stated in Theorem \ref{Baer-exact} and show that in addition to extracting new results from this theorem,  some previous results are obtained from our theorem with simpler proof.  The statements of the following Lemma have routine arguments and so   we leave  their
proofs, but we use them later.\\

\begin{lemma}\label{ann-1} Let $M_R$ be a nonzero module and End$_R(M) = S$. Then the following statements hold.\\
{\rm (i)} $r_M(l_S(r_M(I))) = r_M(I)$ for all  $I\subseteq S$.\\
{\rm (ii)} $l_S(r_M(l_S(N))) = l_S(N)$ for all  $N\subseteq M$.\\
{\rm (iii)} $I^{ece} = I^e$ for all  $I\leq S_S$.\\
{\rm (iv)} $N^{cec} = N^c$ for all  $N\leq M_R$.\\
{\rm (v)} If $I\leq^{\oplus}$$S_S$ (resp. $I\leq^{\oplus}$$_SS$)
then $I^e\leq^{\oplus}$$M_R$
 (resp. r$_M(I)\leq^{\oplus}$$M_R$).\\
{\rm (vi)} If $N\leq^{\oplus}$$M_R$ then
$N^c\leq^{\oplus}S_S$ and l$_S(N)\leq^{\oplus}$$_SS$.\\
{\rm (vii)} $N\in$ Gen$(M_R)$ if and only if $N = I^e$ for some
$I\leq$ $S_S$ if and only if Tr$(M,N) = N$.
\end{lemma}

The following Proposition that may be found
in the literature,  we record that as a  consequence of   Lemma \ref{ann-1}
  for latter uses.\\

\begin{proposition}\label{Baer and dual}{\rm } Let $M_R$ be a nonzero module with End$_R(M) = S$.
Then the following statements hold.\\
{\rm (i)} The  module $M_R$ is Baer if and only if
r$_M(I)$$\leq^{\oplus}M_R$  for every $I\leq$$_SS$. \\
{\rm (ii)} The module $M_R$ is dual Baer if and only if
$IM\leq^{\oplus}$$M_R$  for every $I\leq S_S$.

\end{proposition}

\noindent \emph {\textsf{Proof.}} (i) $(\Rightarrow)$ Let $I\leq$$_SS$ and $N =$ r$_M(I)$. By Lemma \ref{ann-1}(i) we have $N =$
r$_M(l_S(N))$. Now since $M_R$ is Baer $ l_S(N)\leq^{\oplus}$$_SS$. Thus $N\leq^{\oplus} M_R$  by Lemma \ref{ann-1}(v).\\
$(\Leftarrow)$ Let $N\leq M_R$ and $I = l_S(N)$. By our assumption,  r$_M(I)$$\leq^{\oplus}M_R$. Now by parts (vi) and (ii)  of Lemma \ref{ann-1}, we have
$I = l_S(r_M(I)) \leq^{\oplus}$$_SS$.\\
(ii) This is easily obtained  by  parts (iii)-(vi)  of Lemma \ref{ann-1}.\\

\begin{lemma}\label{ann} Let $N\leq M_R$ and End$_R(M) = S$. The following statements are
equivalent.\\
{\rm (i)} $N = $ r$_M(I)$ for some $I\leq$$_SS$.\\
{\rm (ii)} $N =\cap_{ f\in X}\ker f$ for some non-empty $X\subseteq  S$.\\
{\rm (iii)} $M/N\in $ Cog$(M_R)$.\\

\end{lemma}

\noindent \emph {\textsf{Proof.}}   (i)$\Rightarrow$(ii) and (ii)$\Rightarrow$(iii) are clear.\\
(iii)$\Rightarrow$(i). Suppose that $\theta: M/N\rightarrow
M^{\Lambda}$ is an injective $R$-homomorphism for some set
$\Lambda$.  For each $\lambda\in \Lambda$, let $f_{\lambda} =
\pi_{\lambda}\theta$ where $\pi_{\lambda}$ is  the canonical
projection on $M^{\Lambda}$. If $X = \{f_{\lambda} \ | \ \lambda\in \Lambda\}$ then it is easily seen that
 $N =$ r$_M(X)$. Clearly r$_M(X)$ = r$_M(I)$ where $ I = SX $.\\

\begin{theorem}\label{Baer-exact}{\rm }
{\rm (i)} Every exact sequence  $0\rightarrow X\rightarrow
M\rightarrow Y\rightarrow 0$ of $R$-modules
with $Y\in$ Cog$(M)$ splits  if and only if   $M$ is a Baer $R$-module. \\
{\rm (ii)} Every exact sequence  $0\rightarrow X\rightarrow
M\rightarrow Y\rightarrow 0$ of $R$-modules with $X\in$ Gen$(M)$
splits  if and only if   $M$ is a dual Baer $R$-module.
\end{theorem}

\noindent \emph {\textsf{Proof.}} (i) ($\Rightarrow$). To show that  a module $M_R$ is Baer, we use Proposition \ref{Baer and dual}(i).
 Let $N = $ r$_M(I)$ for some $I\leq$$_SS$. By  Lemma \ref{ann}, $M/N\in$ Cog$(M)$. Hence the exact sequence  $0\rightarrow N\stackrel{\imath}\rightarrow
M\stackrel{\pi}\rightarrow M/N\rightarrow 0$ spits by our assumption. This means $N\leq^{\oplus} M$. \\
 ($\Leftarrow$). Consider that if $0\rightarrow
X\stackrel{f}{\rightarrow} M\rightarrow Y\rightarrow 0$ is an
exact sequence of $R$-modules with $Y\in$ Cog$(M_R)$, then the
submodule Im$f$ of $M$ satisfies the equivalent  conditions in the
Lemma \ref{ann}. Thus Im $f$ is a direct summand of $M_R$ by the
Baer condition on $M$. It follows that  the exact sequence splits.\\
(ii) This is dual of (i) and has a similar argument. Just note that  by Lemma \ref{ann-1}(vii),
a submodule $N$ of $M$ is in Gen$(M)$ if and only if $N = IM$ for some $I\leq S_S$.\\

\begin{corollary}\label{baer ann proj } If $R$ is a Baer ring then ann$_R(P)$ is a direct summand of $R_R$ for any projective right $R$-module $P$.
\end{corollary}

\noindent \emph {\textsf{Proof.}} Let $R$ be a Baer ring and $P_R$ be projective. If $I$ = ann$_R(P)$ then $R/I\in $ Cog$(R_R)$. Hence
$I$ is a direct summand of $R_R$ by Theorem
\ref{Baer-exact}.\\

By well known result from Gentile and Levy \cite[Proposition 6.12]{Goodearl-Warfield},  torsionfree  divisible $R$-modules are injective when $R$ is a semiprime right Goldie ring. hence  the following result is a simple corollary of Theorem \ref{Baer-exact}. The part (a) of
below can be  a generalization of \cite[Proposition 4.16]{w.dual baer2015}. An $R$-module $M$ is said to be {\em torsionless} if $M$ is cogenerated by $R$.\\

\begin{corollary}\label{torsion free Baer}{\rm } (a) If $R$ is a semiprime right Goldie (resp. right hereditary Noetherian) ring $R$, every torsionfree injective (resp. injective) module $M_R$ is dual Baer.\\
(b) If $R$ is a ring such that $R^{(n)}_R$ is extending (e.g., $R$ is right self-injective) then  every finitely generated non-singular $R$-module $M_R$ is Baer.\\
(c) If  $R$ is a right hereditary right Noetherian prime ring then every finitely generated  torsionless $R$-module $M$ is Baer.
\end{corollary}

\noindent \emph {\textsf{Proof.}}   Let $0\rightarrow X\rightarrow
M\rightarrow Y\rightarrow 0$ be an exact sequence of $R$-modules
with $X\in$ Gen$(M)$ (resp.  $Y\in$ Cog$(M)$). For (a) (resp. (b)), our assumptions on $R$ imply that $X_R$ is injective (resp. $Y_R$ is projective). For (c), note that $Y$ is also torsionless
 and so by \cite[Proposition 3.4.3]{mc-rob}, $Y_R$ must be projective.  Hence, in any cases,  the exact sequence splits and
the result holds by Theorem \ref{Baer-exact}.\\

\begin{corollary}\label{morita}{\rm } Being Baer (dual Baer) module is a Morita invariant  property.

\end{corollary}

\noindent \emph {\textsf{Proof.}} This is obtained by Theorem
\ref{Baer-exact} and the fact that the category equivalences
preserve exact sequences  and direct (co) products
\cite[Proposition 21.6(3) and (5)]{AF}.\\

\begin{corollary}\label{ch of baer ring}{\rm } A ring $R$ is Baer
if and only if every cyclic torsionless right (left) $R$-module is
projective. Furthermore, every $n$-generated torsionless right(left) $R$-module is projective if and only if M$_n(R)$ is a Baer ring.
\end{corollary}

\noindent \emph {\textsf{Proof.}}  Since a cyclic $R$-module $R/I$ is projective if and only if the exact sequence
$0\rightarrow I\rightarrow R\rightarrow R/I\rightarrow
0$  splits, the result is an application of Theorem
\ref{Baer-exact} for $M = R$. The last statement is now obtained by the fact that the standard Morita equivalent between $R$ and M$_n(R)$ corresponds
 $n$-generated $R$-modules to cyclic M$_n(R)$-modules. \\

\begin{corollary}\label{d.summand of (daul)baer }{\rm }\cite[Theorem  2.17]{rizv1} and \cite[Corollary 2.5]{dual.baer2010}.   A direct summand of a
Baer (resp. dual Baer) module is a Baer (resp. dual Baer) module.
\end{corollary}

\noindent \emph {\textsf{Proof.}} This is obtained by Theorem
\ref{Baer-exact} and the fact that  an exact sequence
$0\rightarrow X\stackrel{f}\rightarrow N\stackrel{g}\rightarrow Y\rightarrow
0$  splits if the exact sequence $0\rightarrow X\oplus
L\stackrel{f\oplus 1_L}\rightarrow N\oplus L\stackrel{\beta}\rightarrow Y\rightarrow 0$ with $\beta(n,l) = g(n)$ splits. In fact,  if there exists $N\oplus
L\stackrel{\alpha}\rightarrow X\oplus L$ such that $\alpha(f\oplus 1_L) =
1_{X\oplus L}$, then we have $hf = 1_X$ where $h:N\rightarrow X$ with
$h(n) = \pi \alpha(n,0)$ and $\pi : X\oplus L\rightarrow X$
 is the natural projection.\\

\begin{lemma}\label{direct-sum-sequence} Let $M = \bigoplus_{\i\in I} M_i$($I$ an index set) such that  Hom$_R(M_i,M_j) = 0$ ($i\neq j$).\\ If  $0\rightarrow X\stackrel{f}\rightarrow M\stackrel{g}\rightarrow Y\rightarrow 0$
is  an exact sequence of $R$-modules. For each $i$, replace $\iota_i(M_i)$ with $M_i$ and let $g_i = g|_{M_i}$,
 $K_i = \ker g_i$, $X_i = f^{-1}(K_i)$ and $f_i = f|_{X_i}$. Then we have:\\
 (a) For each $i$, the sequence  $0\rightarrow X_i\stackrel{f_i}\rightarrow M_i\stackrel{g_i}\rightarrow g(M_i)\rightarrow 0$  is exact.\\
 (b) If  $Y\in$ Cog$(M)$ then for each $i$,   $g(M_i)\in$ Cog$(M_i)$ and  $g(M) = \bigoplus_{\i\in I} g(M_i)$.\\
 (c) If $X\in $ Gen$(M)$ then for each $i$, $X_i\in $ Gen$(M_i)$ and  $X = \bigoplus_{\i\in I} X_i$.\\
 (d) If $N\oplus L = M$ then $N = \bigoplus_{\i\in I} (N\cap M_i)$.
\end{lemma}

\noindent \emph {\textsf{Proof.}} (a) It has a routine argument.\\
(b) Let  $Y\stackrel{\theta}\hookrightarrow \bigoplus_i(M_i)^{\Lambda}$ then by hypothesis, $\theta g(M_i)\subseteq (M_i)^{\Lambda}$ for each $i\in I$.
It follows that $g(M) = \bigoplus_{\i\in I} g(M_i)$ and $g(M_i)\in$ Cog$(M_i)$.\\
(c) Clearly, $\{X_i\}_{i\in I}$ are $R$-linear independent. Let $\bigoplus_i(M_i)^{(\Lambda)}\stackrel{\alpha}\rightarrow X $ be a surjective $R$-homomorphism.
 By hypothesis $f\alpha\psi_i ((M_i)^{(\Lambda_i)} \subseteq M_i$ where $\psi_i:M_i^{(\Lambda_i)}\rightarrow \bigoplus_i(M_i)^{(\Lambda_i)}$ is the natural $R$-monomorphism.
 If $x\in X$, then there are $u_i\in M_i^{(\Lambda_i)}$ $(i = 1, ..., n)$ such that $\alpha(\sum_i u_i) = x$. Since $f\alpha(u_i)\in M_i$, we have $f\alpha(u_i)\in K_i$. Thus
 $x\in \sum_i X_i$. It follows that $X = \sum_i X_i$ and so  $X_i\in $ Gen$(M_i)$.\\
(d) By hypothesis each $M_i$ is a fully invariant submodule of $M_R$. Hence for each $i$ we have $M_i = (N\cap M_i)\oplus (L\cap M_i)$. Thus
$M = [\bigoplus_i (N\cap M_i)] \oplus [\bigoplus_i (L\cap M_i)] = N\oplus L$. It follows that $N = \bigoplus_i (N\cap M_i)$.\\

\begin{corollary}\label{ortognal}{\rm } Let $M = \bigoplus_{\i\in I} M_i$($I$ an index set) such that  Hom$_R(M_i,M_j) = 0$ ($i\neq j$). \\
(a) \cite[Propsitin 3.20]{rizvi2009}. If every $M_i$ is a Baer $R$-module then $M_R$ is Baer.\\
(b) If  every $M_i$ is a dual Baer $R$-module then $M_R$ is dual Baer.
\end{corollary}

\noindent \emph {\textsf{Proof.}}  These are  obtained by Lemma \ref{direct-sum-sequence} and Theorem \ref{Baer-exact}.\\

\noindent  The Theorem \ref{Baer-exact} and Lemma \ref{ann-1}(vii) show that a module $M_R$ is dual Baer
if and only if Tr$(M, X)\leq^{\oplus}M_R$ for every $X\leq M_R$. Clearly, Tr$(M,X)\leq X$. Thus if $M_R$ is dual Baer then a submodule  $X\leq M_R$
contains a non-zero direct summand of $M_R$ if and only if Hom$_R(M, X) \neq 0$. Modules $M_R$ in which  Hom$_R(M,N)\neq 0$
for every $0\neq N\leq M_R$ is called {\em retractable} \cite{khuri2000}. Dually, an $R$-module $M_R$ is called {\em
co-retractable} if Hom$_R(M/N, M)\neq 0$ for every $N < M_R$. In the next section, we shall study the retractable  and co-retractable  conditions for (dual) Baer modules.
 We  introduce the dual notation for Tr$(M,N)$ where $N\leq M$. By  Rej$^{-1}(N)$, we mean  $\pi^{-1}$(Rej$(M/N,M))$ where $\pi : M\rightarrow M/N$ is the
canonical epimorphism. Clearly $N\leq $ Rej$^{-1}(N)$, and $N = $ Rej$^{-1}(N)$
if and only if $M/N\in$ Cog$(M)$. We record the following characterization for Baer (dual Baer) modules and use that in the next section.\\

\begin{proposition}\label{preserv direc-rej-1}{\rm }
 (i) The following conditions are equivalent  for a module $M_R$. \\
 (a) $M_R$ is Baer.\\
 (b) For every nonempty set $\Lambda$ and every
$R$-homomorphism  $f: M \rightarrow M ^{\Lambda}$, the inverse image  $f^{-1}(D)$ is a  direct summand of $M$ where $D$ is a direct summand of $M ^{\Lambda}$.\\
 (c) For every nonempty set $\Lambda$ and every $R$-homomorphism  $f: M \rightarrow M ^{\Lambda}$, $\ker f$ is a direct summand of $M_R$.\\
 (d) For every
$N\leq M_R$,  Rej$^{-1}(N)\leq^{\oplus}M_R$.\\
\indent (ii)  The following conditions are equivalent for a module $M_R$. \\
(a)  $M_R$ is dual Baer.\\
(b) For every nonempty set $\Lambda$, every
$R$-homomorphism $f: M ^{(\Lambda)} \rightarrow M $    preserves  direct summands.\\
(c)  For every nonempty set $\Lambda$, every
$R$-homomorphism $f: M ^{(\Lambda)} \rightarrow M $, Im $f$ is a direct summand of $M_R$.\\
(d)  For every $N\leq
M_R$, Tr$(M,N)\leq^{\oplus}M_R$.
\end{proposition}

\noindent \emph {\textsf{Proof.}} We proof (i). (a)$\Rightarrow$(b). Let $M_R$ be Baer and  $f: M \rightarrow M ^{\Lambda}$ be an $R$-homomorphism for some $\Lambda$.
If $D$ is a direct summand of  $M ^{\Lambda}$ and $N = f^{-1}(D)$ then we have the natural monomorphism $M/N \rightarrow M/D$. It follows that $M/N \in$ Cog$(M)$.
Hence $N$ is a direct summand of $M$ by Theorem \ref{Baer-exact}.  \\
(b)$\Rightarrow$(c). is clear.\\
(c)$\Rightarrow$(a).  Let $\{f_i\}_i \in$ End$_R(M)$. Consider the $R$-homomorphism  $f: M \rightarrow M ^{\Lambda}$ with  $f(m) = \{f_i(m)\}_i$. Then by our assumption
$f^{-1}(0)$ is a direct summand of $M_R$. Thus $\cap_i\ker f_i$ is a direct summand of $M$, proving that $M_R$ is Baer.\\
The equivalence (d)$\Leftrightarrow$ (a), follows by Theorem
\ref{Baer-exact} and the above notes.\\

\section{\bf (Dual) Baer module with retractable or co-retractable condition}

We are now going to study the retractable (co-retractable) condition for Baer and dual Baer modules and then we shall give some applications of our results.\\

\begin{lemma}\label{Rej and tr} Let $M_R$ be a nonzero module.\\
{\rm (i)} If $M_R$ is co-retractable then Rej$^{-1}(N)/N \ll M/N$ for every proper $N\leq M_R$.\\
{\rm (ii)} If $M_R$ is retractable then Tr$(M, N)\leq_{ess}N$ for every nonzero $N\leq M_R$.\\
\end{lemma}

\noindent \emph {\textsf{Proof.}} We only prove (i). Let $N < M_R$
and $K = $ Rej$^{-1}(N)$. Clearly $N\leq K$. Suppose that $K/N +
L/N = M/N$. If $L\neq M$, since $M_R$ is co-retractable, there
exists nonzero homomorphism $g: M/L\rightarrow M$. Consider the
 natural epimorphism $p: M/N\rightarrow M/L$, then $0\neq gp\in$
 Hom$_R(M/N, M)$ and $gp(L/N) = 0$. On the other hand, by the definition of $K$, we have $gp(K/N) = 0$.
 It follows that $gp(M/N) = 0$, a contradiction. Thus $L=M$ and we
 are done.\\

\noindent The equivalences (i)$\Leftrightarrow$ (ii) and (i)$\Leftrightarrow$ (iii) of below are  dual of each others.
 The equivalence (i)$\Leftrightarrow$ (iii) is  appeared in \cite[corollary 2.19]{weak dual}.\\

\begin{theorem}\label{ret and coret}{\rm } The following statements are equivalent for a module $M_R$.\\
{\rm (i)} $M_R$ is semisimple.\\
{\rm (ii)} $M_R$ is  co-retractable and Baer.\\
{\rm (iii)} $M_R$ is  retractable and dual Baer.
\end{theorem}

\noindent \emph {\textsf{Proof.}} We need to show that (ii) or (iii)
$\Rightarrow$ (i). Let  $M_R$ be co-retractable and Baer (the other case is similar). If $N\leq M_R$ and $K$ =  Rej$^{-1}(N)$ then  $K\leq^{\oplus}M_R$ by Proposition \ref{baer and rej}. This shows that
$K/N \leq^{\oplus}M/N$.
On the other hand, $K/N \ll M/N$ by Lemma \ref{Rej and tr}. Hence $N = K$.  Therefore,  every submodule of $M$ is a direct of $M$, as desired.\\

\begin{corollary}\label{qp-dual baer}{\rm } Let $M_R$ be a non-zero quasi-projective module. Then the $R$-module $M$/J$(M)$ is dual Baer if and only if
 it is a semisimple $R$-module.
\end{corollary}

\noindent \emph {\textsf{Proof.}} By \cite[3.4]{Ext} any quasi-projective module with zero Jacobson radical is retractable. Hence the result is obtained by  Theorem \ref{ret and coret}. \\

In \cite{tolooei-vedadi}, rings over which all nonzero modules are
retractable  are studied. Hence over such rings (e.g., commutative
semi-Artinian rings) dual Baer modules are precisely semisimple
modules. \\

\begin{theorem}\label{f.g daul baer}{\rm } Let $R$ be a ring Morita invariant to a  commutative ring.\\
(i) A finitely generated module $R$-module is dual Baer if and only if
it is semisimple.\\
(ii) If $R$ is semi-Artinian then dual Baer $R$-modules are precisely semisimple $R$-modules.
\end{theorem}

\noindent \emph {\textsf{Proof.}} By Corollary \ref{morita}, we can suppose that $R$ is a commutative ring. Thus by \cite[Theorem 2.7]{fin-ret} every finitely generated $R$-module
is retractable. Also if $R$ is semi-Artinian then all nonzero $R$-modules are retractable \cite[Theorem 2.8]{fin-ret}. Thus the result is obtained by Theorem \ref{ret and
coret}.\\

If $X$ and $Y$ are $R$-modules then it is well known that Rej$(X,Y)$ is a fully invariant submodule of $X$ and $X/$Rej$(X,Y)$ lies in Cog$(Y)$. Thus if $M_R$ is a Baer module with no
 non-trivial fully invariant direct summand, then Theorem \ref{Baer-exact} shows that for any $0\neq N\leq M$ we have Hom$_R(M,N) = 0$ or $M\in $ Cog$(N)$. In \cite[Theorem 2.5]{t.dim},
 the {\em triangulating  dimension ($\tau dim(M)$)} was defined for a module $M_R$ as follows:  Sup$\{k \in \Bbb{N} \ | \
M = \oplus_{i=1}^k M_i$ with $M_i\neq 0$ and ${\rm Hom}_R(M_i, M_j) = 0 $ for any $ i < j\}$ and  it was shown that $\tau dim(M_R)<\infty$  if and only if
  $M = \bigoplus_i M_i$ where Hom$_R(M_i, M_j) = 0$ ($i < j$) and each $M_i$ has no non-trivial fully invariant direct summand(i.e.,  $\tau dim(M_i) =1)$ if and only if $M$ has ascending and descending chain conditions  on fully invariant direct summands if and only if
 End$_R(M)$ has  a generalized triangular matrix representation. In \cite[Proposition 2.16]{birken 2000} Baer rings with  a generalized triangular matrix representation are studied. Below we study  Baer modules with
finite   $\tau$dimension and give some applications that are generalizations of earlier results in the literature. Note that
 {\em Noetherian condition $\Rightarrow$ finite uniform dimension  $\Rightarrow$  ascending (descending) chain condition on direct summands $\Rightarrow$ finite   $\tau$dimension}. A modules that  is cogenerated by
each of its nonzero submodule, is called {\em prime} in \cite{bican 1980} and  $*$-{\em prime} in \cite{smithcomp}). It is easy to verify that {\em every prime module  has no non-trivial fully invariant direct summand)}.\\ The following theorem shows that {\em  the study of Baer modules with finite $\tau$dimension  reduces
to the study of such modules when they are prime}. Recall that two $R$-module $X$ and $Y$ are called {\em orthogonal} to each other, if they do not contain nonzero isomorphic submodules. \\

 \begin{lemma}\label{fuldir sumd baer} (a) Let $M_R$ be a nonzero retractable Baer module. Then either $M_R$ is prime or there is a decomposition
$M = N\oplus K$ such that $N$ is a non-trivial fully invariant direct summand of  $M_R$.\\
(b) If $ M = M_1\oplus M_2$ is a Bear $R$-module and $M_i$ is prime $R$-module ($i= 1, 2$) such that Hom$_R(M_1, M_2) = 0$, then
Hom$_R(M_2, M_1) = 0$ and $M_1$ and $M_2$ are orthogonal to each other.
\end{lemma}

\noindent \emph {\textsf{Proof.}} (a) Let $M_R$ is not a prime module. Thus there exists a  non-trivial submodule $X$ of $M$ such that  Rej$(M,X) =: N$ is nonzero.
Clearly $N$ is  a fully invariant  of $M$ such that
$M/N\in$ Cog$(M)$. Since $M_R$ is retractable,  $N\neq M$. Thus $N$ is a non-trivial fully invariant submodule of $M_R$ by  Theorem \ref{Baer-exact}.\\
(b) If $f : M_2 \rightarrow M_1$ is nonzero then by the Baer condition on $M$, we must  have  $M_2 \simeq \ker f\oplus $Im $f$. Now since $M_1$ is prime, it lies in Cog(Im $f$). Hence
Hom$_R(M_1, M_2) \neq 0$, contradiction.\\

 \begin{theorem}\label{baer t.dim}{\rm }
Let $M_R$ be a non-zero Baer module. Then $M = \bigoplus_{i=1}^n M_i$  such that
each $M_i$  is a prime $R$-module and   $M_i's$ are mutually orthogonal with
Hom$_R(M_i, M_j) = 0$ ($i \neq j$) if and only if
 $\tau dim(M_R)<\infty$  and  every direct summand of $M_R$ is retractable.
\end{theorem}

\noindent \emph {\textsf{Proof.}} ($\Leftarrow$).   This obtained by  \cite[Theorem 2.5]{t.dim} and Lemma \ref{fuldir sumd baer}.\\
 ($\Rightarrow$). Let $0\neq K\leq N\leq M = N\oplus N'$. We shall show that  Hom$_R(N,K)\neq 0$. By induction, we can suppose that $N\cap M_i \neq 0$ for all $i$. Also by Lemma \ref{direct-sum-sequence}(d),
we have $N =  \bigoplus_{i=1}^n (N\cap M_i)$. On the other hand, we may suppose that $K\cap M_j \neq 0$ for some $j$. Thus the prime condition on $M_j$ implies that
Hom$_R(N\cap M_j, K) \neq 0$, proving that  Hom$_R(N,K)\neq 0$.\\

Following \cite{smithcomp})  a module $M_R$ is called {\em compressible} if  $M$ can be embedded in every non-zero its submodule.
$R$-modules $X$ and $Y$ are said to be {\em sub-isomorphic} if $X$ can be embedded in $Y$ and vise versa.\\

\begin{proposition}\label{baer ret}{\rm }
Let $M_R$ be a non-zero Baer module with ascending (descending) chain condition on direct summands. Then $M = \bigoplus_{i=1}^n M_i$  such that
Hom$_R(M_i, M_j) = 0$ ($i\neq j$) and each $M_i$  is a finite direct sum of indecomposable compressible $R$-modules that are  mutually sub-isomorphic
 if and only if  every direct summand of $M_R$ is retractable.

\end{proposition}

\noindent \emph {\textsf{Proof.}}  ($\Rightarrow$). It is easy to show  that each $M_i$ is a prime .  Hence every direct summand of $M_R$ is retractable by Theorem \ref{baer t.dim}.\\
($\Leftarrow$). By \cite[Proposition 6.59]{Lam2}, $\tau dim(M_R)$ is finite. Hence by Theorem \ref{baer t.dim}, we may suppose that $M_R$ is a prime module. Thus Hom$_R(X,Y) \neq 0$ for every non zero submodules $X$ and $Y$ of $M_R$.
 On the other hand, by \cite[proposition 10.14]{AF},  $M = \bigoplus_{i=1}^m P_i$ is a finite direct sum of indecomposable  submodules.
 Since now $P_i\oplus P_j$ is Baer for all $i, j$,  every $f: P_i\rightarrow P_j$ is  one to one. It follows that $M$ must be a finite direct sum of indecomposable compressible $R$-modules that are  mutually sub-isomorphic.\\

As we stated in the proof of Theorem \ref{f.g daul baer}, if $R$ is commutative and $M_R$ is finitely generated then every direct summand of $M_R$ is retractable.
 Thus the following result can be  a generalization of \cite[Proposition 2.19]{rizv1}.\\

\begin{theorem}\label{comm-domain fg}{\rm } Let  $R$ be a ring Morita invariant to a  commutative domain.
Then every finitely generated Baer $R$-module $M$ with a finite $\tau$dimension  is either torsion or    a finite direct sum of uniform right ideals of $R$.
\end{theorem}

\noindent \emph {\textsf{Proof.}} First note that $M_R$  is singular if and only
if there exists an exact sequence $0\rightarrow  B \rightarrow A \rightarrow M \rightarrow 0$ such that the map $B\rightarrow A$ is
an essential monomorphism. Since (essential) monomorphisms and co-kernels are preserved under Morita equivalences,  we may suppose by Corollary \ref{morita} that $R$ is a commutative domain. Thus every prime $R$-module is torsion or torsionless (note that if $M_R$ is not torsion
then $M$ contains  isomorphically a nonzero (right) ideal of $R$). Also since $R$ is assumed to be a commutative domain, the uniform dimension of $R$ is finite. It follows that
 every torsionfree $R$-module is faithful with  finite  uniform dimension.  On the other hand, by \cite[Theorem 3.14]{smith-vedadi weakg}, if $M_R$ is  finitely generated and faithful then Hom$_R(M, X)\neq 0$ for every nonzero $X_R$. Now  we can apply Theorem \ref{baer t.dim}  to deduce that $M_R$  is  either torsion or faithful torsionfree $R$-module  with  finite  uniform dimension. If
  $M_R$  is  faithful torsionfree $R$-module  with  finite  uniform dimension then Proposition \ref{baer ret} shows that $M$ is a finite direct sum of uniform right ideals of $R$.  \\

\begin{corollary}\label{comm-her-domain fg}  Let  $R$ be a ring Morita invariant to a  commutative hereditary Noetherian domain and $M_R$ is  finitely generated.
Then $M_R$ is a Baer $R$-module $M$ if and only if it is semisimple  or   projective.

\end{corollary}

\noindent \emph {\textsf{Proof.}}  The sufficiency follows from Corollary \ref{torsion free Baer}. For the necessity,
 note that every singular $R$-module has nonzero socle by \cite[Proposition 5.4.5]{mc-rob}. Hence the result is obtained
by Proposition \ref{baer ret} and Theorem \ref{comm-domain fg}.\\

 We now consider  the dual of Theorem \ref{baer t.dim}. Following   \cite {wiz coprime}, a  module $M$ is called  {\em  co-prime}  whenever $M/N$ generates $M$ for any proper submodule $N< M_R$. Furthermore,
  we say that $M_R$ is {\em co-compressible} if $M_R$ is a homomorphic image of every non-zero its factor. Clearly,   co-compressible modules are    co-prime.  For the dual notion of
  sub-isomorphic, by $X\stackrel{epi}{\simeq}Y$ the $R$-modules $X$ and $Y$ are called {\em epi-invariant}  if $X$ is a homomorphic image of $Y$ and vise versa.\\

 \begin{lemma}\label{fuldir sumd dualbaer} (a) Let $M_R$ be a nonzero co-retractable dual Baer module. Then either $M_R$ is co-prime or there is a decomposition
$M = N\oplus K$ such that $N$ is a non-trivial fully invariant direct summand of  $M_R$.\\
(b) If $ M = M_1\oplus M_2$ is a dual Bear $R$-module and $M_i$ is co-prime $R$-module ($i= 1, 2$) such that Hom$_R(M_1, M_2) = 0$, then
Hom$_R(M_2, M_1) = 0$. \\
(c) Let $M = \bigoplus_{i\in I} C_i$ such that for any $i\in I$, the $R$-module $C_i$ has no non-trivial fully invariant direct summand. If $V$ is a non-zero fully invariant direct summand of $M_R$
then there is non-empty subset $J$ of $I$ such that  $ V = \bigoplus_{j\in J} C_j$ and Hom$_R(V, C_i) = 0$ for any $i\in I\setminus J$.
\end{lemma}

\noindent \emph {\textsf{Proof.}} (a) and (b) are  dual of  Lemma \ref{fuldir sumd baer} and  obtained   from the dual arguments.\\
(c) Since $V$ is fully invariant submodule of $M_R$,  we must have $V = \bigoplus_{i\in I} V\cap C_i$, also if $M = V\oplus V'$ then Hom$_R(V,V') = 0$.
On the other hand, for every $i\in I$, it is easy to verify that
 $V\cap C_i$ is a fully invariant direct summand of $C_i$.  Now let $J = \{ i\in I \ | \ V\cap C_i \neq 0 \}$.\\

 \begin{theorem}\label{dual baer t.dim}{\rm }
Let $M_R$ be a non-zero dual Baer module. Consider the following conditions.\\
(a) Every direct summand of $M_R$ is co-retractable.\\
(b) $M = \bigoplus_{i\in I} M_i$ is a  direct sum of co-prime $R$-modules with Hom$_R(M_i, M_j) = 0$ ($i \neq j$) such that
each $M_i$ is a  direct sum of  indecomposable co-compressible modules   that are  mutually epi-invariant.\\
Then (a) implies (b) and the converse is true if $I$ is a finite set.
\end{theorem}

\noindent \emph {\textsf{Proof.}} (a)$\Rightarrow$(b). By   \cite[corollary 2.6]{dual.baer2010}, $M = \bigoplus_{ i\in I} C_i$ is a direct sum of indecomposable modules.
Also every indecomposable  co-retractable dual Baer module $X$ is a co-compressible module. To see this suppose that there is a proper submodule $Y$ of $X$. Since $X$ is assumed to be co-retractable,
 there exits  non-zero $f : X/Y \rightarrow X$. Now Im$(f)$ lies in Gen$(X)$ and so must be a direct summand of $X$ by the dual Baer condition. It follows that $f$ is surjective,
 proving that $X$ is   co-compressible. Now for any $i\in I$, let $M_i = \bigoplus\{C_j \ | \  C_i\stackrel{epi}{\simeq}C_j\}$. By Lemma \ref{fuldir sumd dualbaer}(c), $M_i$ has no non-zero
  fully invariant direct summand and hance by part (a) of Lemma \ref{fuldir sumd dualbaer},  $M_i$ must be a co-prime $R$-module. On the other hand, if  Hom$_R(C_i, C_j) \neq 0$ ($i \neq j$) then
   $ C_i\stackrel{epi}{\simeq}C_j$. It follows that  Hom$_R(M_i, M_j) = 0$ ($i \neq j$). \\
   For the converse, let $I$ be a finite set with $|I| = n$,  $N\oplus L = M$ and $K < N_R$. We shall show that Hom$_R(N/K, N) \neq 0$.
   Let $V= K\oplus L$. Then $N/K \simeq M/V$ and we will show that   Hom$_R(M/V, N) \neq 0$. In fact, it is enough to show that there
   exists $i\in I$ such that $M_i\not\subseteq V$ and  $M_i\in$ Gen($M/V$). Because by hypothesis $M_i$ is a fully invariant submodule of $M_R$
     and so we  have $M_i = (N\cap M_i)\oplus(L\cap M_i)$ with $N\cap M_i \neq 0$.  Now let  $ W = M_1 + V$. If $W = M$ then $M_1\not\subseteq V$ and $M/V\simeq M_1/(V\cap M_1)$.
     Hence the co-prime condition on $M_1$ implies that $M_1\in$ Gen($M/V$), as desired. If $W\neq M$,  let
    $J = \{i\in I \ | \ M_i\subseteq W\}$ and $U = \bigoplus_{j\in J} M_j$.  Then $W/U$ is a proper submodule of $M/U\simeq \bigoplus_{i\not\in J} M_i$ which is
   co-retractable  by the induction assumption. Thus    Hom$_R(M/W, M_i) \neq 0$ for some $i\not\in J$.  It follows Hom$_R(M/W, C) \neq 0$ such that  $C$ is an indecomposable  module  appeared in the decomposition of $M_i$.
    Now the dual Baer condition on $M_R$ implies that
   $C$ is a homomorphic image of $M/W$. Since   indecomposable  modules  in the decomposition of $M_i$  are  mutually epi-invariant, we have $M_i\in$ Gen($M/W) \subseteq$ Gen$(M/V)$.
   The proof is now completed.\\

\begin{proposition}\label{d-bear co-ret}{\rm}
Let $M$ be a non-zero dual Baer $R$-module  such that every direct summand of $M_R$ is co-retractable. Then
$M = M_1\oplus M_2$ where $M_1$ is singular and $M_2$ is a nonsingular semisimple module.
\end{proposition}

\noindent \emph {\textsf{Proof.}} If $X$ is co-prime and $N$ is a proper essential submodule of $X$, then  $X/N$ and hence $X$ must be singular.
 Therefore every co-prime module $X$ is either singular or semisimple (as semisimple modules have no proper essential submodules). The proof is now completed by Theorem \ref{dual baer t.dim}.\\

\begin{theorem}\label{perfect- (dual)baer}{\rm } Let $R$ be a ring Morita invariant to a right duo perfect ring, then the following condition are equivalent.\\
(i) $M_R$ is Baer.\\
(ii) $M_R$ is dual Baer.\\
(iii) $M_R$ is semisimple.

\end{theorem}

\noindent \emph {\textsf{Proof.}} By Corollary \ref{morita}, we  suppose that $R$ is a right duo  and a perfect ring. Then it is well known that every non-zero $R$-module has a non-zero maximal submodule. It follows that
co-compressible $R$-modules are simple. Therefore, by Theorems \ref{ret and coret} and \ref{dual baer t.dim}, it is enough to show that every non-zero $R$-module is co-retractable.
Now let $N$ be a proper submodule of a non-zero module $M_R$. Since $R$ is right perfect then the non-zero module $M/N$ has a maximal submodule $K/N$. On the other hand,
by \cite[Theorem 2.14]{fully kash} $M_R$ is a Kasch module. Hence the simple module $M/K$ can be embedded in $M_R$. It follows that Hom$_R(M/N, M)$ is nonzero, as desired.\\

\section{\bf Further applications}

As another application of the Theorem \ref{Baer-exact}, we conclude the paper with a result to show that if Hom$_{\Bbb Z}(M,{\Bbb Q}/{\Bbb Z})$ is Baer as a left $R$-module then
$M_R$  has a condition close to the dual Baer condition. A submodule $N$ of a module $M_R$ is called {\em pure} if for any
left $R$-module $A$, the homomorphism $i\otimes 1_A$ is one to one
where $i: N\rightarrow M$ is the inclusion map and $1_A:
A\rightarrow A$ is the identity map. The exact sequence
$0\rightarrow X\stackrel{f}{\rightarrow} M\rightarrow Y\rightarrow
0$ is then called {\em pure exact} if Im$f$ is a pure submodule of
$M$. Clearly, every direct summand of $M_R$ is a pure submodule of
$M$. Hence, in view of Proposition \ref{baer and rej}, we may
consider the condition weaker than the dual Baer condition for a module $M$: {\em  Tr$(M,X)$  is a pure submodule $M_R$
for any $X\leq M_R$}.   For any right $R$-module $M$, the character left
$R$-module Hom$_{\Bbb Z}(M,{\Bbb Q}/{\Bbb Z})$ is denoted by
$M^+$.  First we recall  some facts on purity. \\

\begin{proposition}\label{regular-1}{\rm }
Let $M$ be a non-zero $R$-module.\\
{\rm (i)} The exact sequence $0 \rightarrow N\rightarrow M$ of
$R$-modules is pure exact if and only if the sequence
$M^+\rightarrow N^+\rightarrow 0$ splits. \\
{\rm (ii)} Let $M_R$ be flat. An exact sequence $0\rightarrow
X\rightarrow M\rightarrow Y\rightarrow 0$ of $R$-modules  is pure
if and only if $Y_R$ is flat.
\end{proposition}

\noindent \emph {\textsf{Proof.}} (i) This follows from  \cite[
Proposition 5.3.8]{enoch} \\
(ii) It is true by  \cite[Corollary 4.86]{Lam2}.\\

\begin{theorem}\label{M+ Baer}{\rm }
Let $M_R$ be a non-zero module with End$_R(M) = S$. If  the left
$R$-module $M^+$ is Baer then Tr$(M,X)$  is a pure submodule $M_R$
for any $X\leq M_R$. If further, $M_R$ is flat then $M/N$
is a flat $R$-module for any $N\in$ Gen$(M)$.
\end{theorem}

\noindent \emph {\textsf{Proof.}} Suppose that the left $R$-module
$M^+$ is Baer. Let $X\leq M_R$ and Tr$(M,X) = N$. Then
$N\in$ Gen$(M)$. Consider the exact sequence $0\rightarrow
N\rightarrow M\rightarrow M/N\rightarrow 0$ in Mod-$R$. Note that
since $N\in$ Gen$(M)$, we have $N^+\in Cog(M^+)$. Thus  we obtain
the exact sequence $0\rightarrow (M/N)^+\rightarrow M^+\rightarrow
N^+\rightarrow 0$ in $R$-Mod with $N^+\in Cog(M^+)$ by \cite[
Proposition 4.8]{Lam2}. Now by the Baer condition  on $M^+$ and Theorem
\ref{Baer-exact}, the sequence $M^+\rightarrow N^+\rightarrow 0$
splits. Hence $N$ is a pure submodule of $M_R$ by Proposition
\ref{regular-1}(i). The last statement is now obtained by
Proposition \ref{regular-1}(ii). \\

\begin{corollary}\label{comp-pure}{\rm } Let $M$ be a non-zero $R$-module. (i) If $M_R$ is  self-generator and the left
$R$-module $M^+$ is Baer then all submodules of $M_R$  are  pure.\\
(ii)  If $M_R$ is  a generator  for Mod-$R$, then  the left
$R$-module $M^+$ is Baer if and only if it is a semisimple left $R$-module.
\end{corollary}

\noindent \emph {\textsf{Proof.}} (i) This is an immediate corollary of Theorem \ref{M+ Baer}.\\
(ii) Just note that for every left $R$-module $L$, we have $L$ can be embedded in the left $R$-module $L^{++}$. Hence. if
   $M_R$ is a  generator  for Mod-$R$, then $M^+$ is a co-generator for $R$-Mod. Thus the result is obtained by Theorem \ref{Baer-exact}.\\

A ring $R$ is called {\em von Neumann regular} if for every $a\in R$ there is $b\in R$ such that $a = aba$. It is well known that $R$ is a von Neumann regular ring if and only if
 all cyclic (left) right $R$-module are flat if and only if every right (left) ideal is a pure in $R_R$ ($_RR$).\\

\begin{corollary}\label{Baer-Regular}{\rm } If  the left
$R$-module $R^+$ is Baer then $R$ is a von Neumann regular ring.
\end{corollary}

\noindent \emph {\textsf{Proof.}} It is obtained by Corollary \ref{comp-pure} and the above notes.\\

\vskip 0.4 true cm



\end{document}